\documentclass[12pt, draft]{amsart}

\usepackage{amsmath}
\usepackage{amssymb}
\usepackage{amsthm}
\usepackage{bm}
\numberwithin{equation}{section}
\newtheorem{theorem}{Theorem}[section]
\newtheorem{lemma}[theorem]{Lemma}

\newtheorem{corollary}[theorem]{Corollary}

\newtheorem{thm}{Theorem}

\begin{document}
\title[The generalized strong recurrence]{The generalized strong recurrence for non-zero rational parameters}
\author[T.~Nakamura]{Takashi Nakamura}
\address{Department of Mathematics Faculty of Science and Technology Tokyo University of Science Noda, CHIBA 278-8510 JAPAN}
\email{nakamura\_takashi@ma.noda.tus.ac.jp}

\maketitle
\begin{abstract}
The strong recurrence is equivalent to the Riemann hypothesis. On the other hand, the generalized strong recurrence holds for any irrational number. In this paper, we show the generalized strong recurrence for all non-zero rational numbers. Moreover, we prove that the generalized strong recurrence in the region of absolute convergence holds for any real number.
\end{abstract}

\section{Introduction and main result}
The value-distribution of the Riemann zeta function $\zeta (s)$ has been investigated by many mathematicians (see for example \cite{Edwards} and \cite{Tit}). In 1975, S.~M.~Voronin \cite{Voronin} established the universality theorem. Roughly speaking, this theorem implies that any non-vanishing analytic function can be uniformly approximated by the Riemann zeta function $\zeta (s)$. 

To state it, we need some notation.  By ${\rm{meas}} \{A\}$ we denote the Lebesgue measure of the set $A$, and, for $T>0$, we use the notation $\nu_T \{ \ldots \} := T^{-1} {\rm{meas}} \{ \tau \in [0,T] : \, \ldots \}$ where in place of the ellipsis, some condition satisfied by $\tau$ is to be written. Let $D:= \{s \in {\mathbb{C}} : 1/2 < \Re (s) < 1\}$ and $K$ be a compact subset of the critical strip $D$ with connected complement. The strongest version of Voronin's universality theorem is as follows. 
\begin{thm}{\rm{(see \cite[Theorem 6.5.2]{Lali}). }}
Let $f(s)$ be a non-vanishing continuous function on $K$ which is analytic in the interior of $K$. Then for every $\varepsilon > 0$, it holds that
$$
\liminf_{T \rightarrow \infty} \nu_T \Bigl\{  \sup_{s \in K} \bigl| \zeta (s+ i\tau ) - f(s) \bigr| < \varepsilon \Bigr\} > 0 .
$$
\label{th:Vo}
\end{thm}

B.~Bagchi proved that the Riemann hypothesis is true if and only if the Riemann zeta-function  can be approximated by itself in the sense of universality. This property is called strong recurrence  (see \cite[Theorem 8.3]{Steuding2}).  
\begin{thm}{\rm{(see \cite[Theorem 3.7]{BagchiZ}).}}
The Riemann hypothesis is true if and only if, for any $K$ and for any $\varepsilon> 0$, 
\begin{equation*}
\liminf_{T \rightarrow \infty} \nu_T \Bigl\{ \sup_{s \in K} \bigl| \zeta (s+i\tau) - \zeta (s) \bigr| < \varepsilon\Bigr\} > 0 .
\end{equation*}
\label{th:bage}
\end{thm}

Meanwhile, the M\"{o}bius $\mu$-function is defined by $\mu (1) =1$, $\mu (n) = 0$ if $n$ has a quadratic divisor $\ne 1$, and $\mu (n) =(-1)^r$ if $n$ is the product of $r$ distinct primes. The Riemann hypothesis is equivalent to the estimate $M(x) := \sum_{m \le x} \mu (n) \ll x^{1/2+\varepsilon}$ (see \cite[Section 14.25]{Tit}). 

 Denjoy \cite{Denjoy} argued as follows (see also \cite[Section 3.3]{Steuding2}). Assume that $\{ X_n \}$ is a sequence of random variable with distribution ${\rm{P}} (X_n =1) = {\rm{P}} (X_n = -1) =1/2$. Define $S_0 =0$ and $S_n = \sum_{m=1}^n X_m$. By central limit theorem, we obtain $\lim_{n \to \infty} {\rm{P}} (|S_n| \ll n^{1/2+\varepsilon}) =1$. In words of Edwards \cite[Section 12.3]{Edwards}: `Thus these probabilistic assumptions about the values of $\mu (n)$ lead to the conclusion, ludicrous as it seems, that $M(x) = O(x^{1/2+\varepsilon})$ with probability one and hence that the Riemann hypothesis is true with probability one!'. 

Inspired by Theorem \ref{th:bage} and Denjoy's probabilistic argument, the author proved the following generalized strong recurrence. 
\begin{thm}{\rm{(see \cite[Corollary 1.4]{Nakamura1}).}}
For almost all $d \in {\mathbb{R}}$ and for any $\varepsilon > 0$ and $K$, 
\begin{equation}
\liminf_{T \rightarrow \infty} \nu_T \Bigl\{ \sup_{s \in K} \bigl| \zeta (s+ i \tau) - \zeta (s+ id \tau) \bigr| < \varepsilon \Bigr\} > 0 .
\label{eq:3}
\end{equation}
\label{th:na1}
\end{thm}
The author also proved the generalized strong recurrence for every algebraic irrational number in \cite[Corollary 1.2]{Nakamura1}. Afterwards, Pa\'{n}kowski showed the generalized strong recurrence for any irrational number.

\begin{thm}{\rm{(see  \cite[Theorem 1.1]{Pan}).}}
For any irrational $d \in {\mathbb{R}}$ and for any $\varepsilon > 0$ and $K$, 
\begin{equation}
\liminf_{T \rightarrow \infty} \nu_T \Bigl\{ \sup_{s \in K} \bigl| \zeta (s+ i \tau) - \zeta (s+ id \tau) \bigr| < \varepsilon \Bigr\} > 0 .
\label{eq:4}
\end{equation}
\label{th:pa1}
\end{thm}

In the present paper, we will show that the generalized strong recurrence holds for all non-zero rational numbers. The keys of the proofs of Theorems \ref{th:na1} and \ref{th:pa1} are the lemmas similar to Lemma \ref{lem:pan}. Obviously $\{ \log p_n \} \cup \{ \log p_n^d \}$ is not linearly independent over ${\mathbb{Q}}$ for any rational number $d$. Therefore the proof of the following theorem is completely different.
\begin{theorem}
Let $k$ and $j$ be coprime integers. For any $\varepsilon > 0$ and $K$, 
\begin{equation}
\liminf_{T \rightarrow \infty} \nu_T \Bigl\{ \sup_{s \in {\mathcal{K}}} \bigl| \zeta (s+ i j \tau) - \zeta (s+ ik \tau) \bigr| < \varepsilon \Bigr\} > 0 .
\label{mtheq:1}
\end{equation}
\label{mth:m1}
\end{theorem}

Hence by putting $\tau' = j \tau$ in (\ref{mtheq:1}) and using Theorem \ref{th:pa1}, we have;

\begin{corollary}
For any $ 0 \ne d \in {\mathbb{R}}$ and for any $\varepsilon > 0$ and $K$, 
\begin{equation}
\liminf_{T \rightarrow \infty} \nu_T \Bigl\{ \sup_{s \in K} \bigl| \zeta (s+ i \tau) - \zeta (s+ id \tau) \bigr| < \varepsilon \Bigr\} > 0 .
\label{eq:2}
\end{equation}
\label{cor:pa1}
\end{corollary}

 Note that we can prove the above results for some large class of zeta functions which have the Euler product (see for example \cite[Section 2]{Nakamura2}, \cite[Section 2]{Pan} or \cite[Section 2.2]{Steuding2}). 
 
 This paper is divided into 4 Sections. In Section 2 we show the limit theorem to prove Theorem \ref{mth:m1}. In Section 3, we prove Theorem \ref{mth:m1}. In Section 4, we will show that the generalized strong recurrence in the region of absolute convergence holds for any $d \in {\mathbb{R}}$. 

\section{Limit theorem}
To prove Theorem \ref{mth:m1}, we show the Limit theorem \ref{th:lim1}. Denote by $H(D)$ the space of analytic functions on $D$ equipped with the topology of uniform convergence on compacta. We denote by $\gamma$ the unit circle on ${\mathbb{C}}$, and let 
$$
\Omega := \prod_{p} \gamma (p) ,
$$
where $\gamma (p) = \gamma$ for all primes $p$. With the product topology and pointwise multiplication, the infinite dimensional torus $\Omega$ is a compact topological Abelian group. Denoting by $m_{H}$ the probability Haar measure on $(\Omega , {\mathfrak{B}}(\Omega))$, we obtain a probability space $(\Omega , {\mathfrak{B}}(\Omega) , m_{H})$. Let $\omega (p)$ stand for the projection of $\omega \in \Omega$ to the coordinate space $\gamma (p)$. Further, let $\omega (1) =1$ and $
\omega (n) := \prod_p \omega (p)^{v (n;p)}$, where $n \in {\mathbb{N}}$ and $v (n;p)$ is the exponent of the prime $p$ in the prime factorization of $n$. For $s \in D$ and $\omega \in \Omega$, we define
\begin{equation}
\zeta (s, \omega) := \sum_{n=1}^{\infty} \frac{\omega(n)}{n^s} .
\label{eq:zetao}
\end{equation}
This is equal to $\zeta (s)$ for $\sigma >1$ and $\omega \equiv 1$. Let $k$ and $j$ be coprime integers as introduced in Section 1. 

\begin{lemma}
The function $\zeta (s, \omega^j) - \zeta (s, \omega^k)$ is an $H(D)$-valued random element on the probability space $(\Omega , {\mathfrak{B}}(\Omega) , m_{H})$. 
\label{lem:1}
\end{lemma}
\begin{proof}
By \cite[Lemma 4.1]{Steuding2}, the function $\zeta (s, \omega^j)$ converges uniformly on compact subsets of $D$ for almost all $\omega \in \Omega$. Since the finite union of null sets is also a null set, the functions $\zeta (s, \omega^j)$ and $\zeta (s, \omega^k)$ converge uniformly for almost all $\omega \in \Omega$. Hence $\zeta (s, \omega^j) - \zeta (s, \omega^k)$ converges uniformly on compact subsets of $D$ almost surely. Therefore the assertion of the lemma follows. 
\end{proof}

By the decomposition of $\omega$, we may rewrite the series (\ref{eq:zetao}) as
\begin{equation}
\zeta (s, \omega) = \prod_p \biggl( 1 - \frac{\omega (p)}{\smash{p^s}} \biggr)^{-1} =
\prod_p \biggl( 1 + \sum_{m=1}^{\infty} \frac{\omega^m (p)}{\smash{p^{ms}}} \biggr) .
\label{eq:lem2}
\end{equation}
In the half-plane $\sigma >1$ both the series (\ref{eq:zetao}) and the product (\ref{eq:lem2}) converges absolutely. The next lemma is proved by \cite[Lemma 4.2]{Steuding2} and the manner similar to the proof of Lemma \ref{lem:1}. 
\begin{lemma}
For almost all $\omega \in \Omega$,  $\zeta (s, \omega^j) - \zeta (s, \omega^k)$ written by the product (\ref{eq:lem2}) converges uniformly on compact subsets of $D$. 
\label{lem:2}
\end{lemma}

In order to establish limit theorems for Dirichlet polynomials, let
$$
f_N (s) := \sum_{n=1}^N \frac{1}{n^s}, \qquad f_N (s, g) = \sum_{n=1}^N \frac{g(n)}{n^s},
$$
where $g$ is a unimodular, completely multiplicative arithmetic function; in particular, $f_N (s,1) = f_N (s)$, where $1 \colon {\mathbb{N}} \to {\mathbb{C}}$ is the arithmetic function constant $1$. 
\begin{lemma}
Let $P_{T,N}^g (A) := \nu_T \{ f_N (s+i\tau, g^k) - f_N (s+i\tau, g^j) \in A\}$, where $A \in {\mathfrak{B}}(H(D))$. Then there exists a measure $P_N^g$ on $(H(D), {\mathfrak{B}}(H(D)))$ such that $P_{T,N}^g$ converges weakly to $P_N^g$ as $T \to \infty$. 
\label{lem:3}
\end{lemma}
\begin{proof}
Define $\Omega_r := \prod_{m=1}^{r} \gamma(p_m)$ and $h (x_1, \ldots, x_r \,; g^j)$ $\colon \Omega_r \to H(D)$ by
$$
h (x_1, \ldots, x_r \,; g^j) := \sum_{n=1}^N \frac{g^j(n)}{n^s} \prod_{m=1}^r x_m^{-j v(n;p)} .
$$
By using $h (x_1, \ldots, x_r \,; g^j)- h (x_1, \ldots, x_r \,; g^k)$ instead of $h (x_1, \ldots, x_r)$ defined in \cite[p.~68]{Steuding2} and modifying the proof of \cite[Lemma 4.4]{Steuding2}, we obtain this lemma.
\end{proof}

By $h (x_1, \ldots, x_r \,; g^j)- h (x_1, \ldots, x_r \,; g^k)$ and the way similar to the proof of \cite[Lemma 4.5]{Steuding2}, we obtain the following lemma.
\begin{lemma}
Both probability measures $P_{T,N}^g$ and $P_{T,N}^1$ converge weakly to the same measure as $T \to 0$. 
\end{lemma}

Let $s \in D$, $N \in {\mathbb{N}}$ and $\sigma_1 > 1/2$, and put
$$
\zeta_N (s) := \sum_{n=1}^{\infty} \frac{1}{n^s} \exp \bigl( - (n/N)^{\sigma_1} \bigr) ,
$$
where the series converges absolutely for $\sigma > \sigma_1$. In view of \cite[Lemma 4.8]{Steuding2}, the function $\zeta_N (s)$ can approximate $\zeta (s)$ in an appropriate mean. Hence we have the following lemma by using the triangle inequality.
\begin{lemma}
Let ${\mathtt{K}}$ be a compact subset of $D$. Then
$$
\lim_{N \to \infty} \! \limsup_{T \to \infty} \int_0^T \max_{s \in \mathtt{K}}
\bigl| \zeta_N (s+ij\tau) - \zeta_N (s+ik\tau) - \zeta (s+ij\tau) + \zeta (s+ik\tau) \bigr| d\tau =0. 
$$
\label{lem:mean}
\end{lemma}

Next we define the functions for $s \in D$, $\sigma_1 > 1/2$ and $\omega \in \Omega$,
$$
\zeta_{N,M} (s) := \sum_{n=1}^M \frac{1}{n^s} \exp \bigl( - (n/N)^{\sigma_1} \bigr) , \qquad 
\zeta_N (s, \omega) := \sum_{n=1}^{\infty} \frac{\omega (n)}{n^s} \exp \bigl( - (n/N)^{\sigma_1} \bigr) .
$$
\begin{lemma}
Define the probability measures
\begin{equation*}
\begin{split}
&P_{T,N} (A) := \nu_T \{ \zeta_N (s+ij\tau) - \zeta_N (s+ik\tau) \in A \} , \\
&Q_{T,N} (A) := \nu_T \{ \zeta_N (s+ij\tau, \omega^j) - \zeta_N (s+ik\tau, \omega^k) \in A \} ,
\end{split}
\end{equation*}
for $A \in {\mathfrak{B}}(H(D))$. Then there exists a measure $P_N'$ on $(H(D), {\mathfrak{B}}(H(D)))$ such that the measures $P_{T,N}$ and $Q_{T,N}$ converge weakly to $P_N'$ as $T \to \infty$. 
\end{lemma}
\begin{proof}
Let $\theta$ be a random variable uniformly distributed on $[0,1]$, defined on some probability space $({\mathbb{R}}, {\mathfrak{B}}({\mathbb{R}}), P^*)$. Define
$$
{\mathcal{Z}}_{N,M} (s+i\theta T) := \zeta_{N,M} (s+ij\theta T) - \zeta_{N,M} (s+ik\theta T).
$$
By considering ${\mathcal{Z}}_{N,M} (s+i\theta T)$ instead of \cite[(4.23)]{Steuding2}, and modifying the proof of \cite[Lemma 4.9]{Steuding2}, we obtain this lemma.
\end{proof}

Hence we have the following limit theorem by using Lemma \ref{lem:mean} and the method similar to the proof of \cite[Theorem 4.3 and Lemma 4.11]{Steuding2}. 
\begin{theorem}
Define the probability measures
\begin{equation*}
\begin{split}
&P_T (A) := \nu_T \{ \zeta (s+ij\tau) - \zeta (s+ik\tau) \in A \} , \\
&P (A) := m_H \{ \zeta (s, \omega^j) - \zeta (s, \omega^k) \in A \} 
\end{split}
\end{equation*}
for $A \in {\mathfrak{B}}(H(D))$. Then the measure $P_T$ converges weakly to $P$ as $T \to \infty$. 
\label{th:lim1}
\end{theorem}

\section{Proof of Theorem \ref{mth:m1}}
Recall that the minimal closed set $S_{\rm{P}} \subseteq H (D)$ such that ${\rm{P}} (S_{\rm{P}})=1$ is called the support of ${\rm{P}}$. The set $S_{\rm{P}}$ consists of all $x \in H (D)$ such that for every neighborhood $V$ of $x$ the inequality ${\rm{P}} (V) > 0$ is satisfied. The support of the distribution of the random element $X$ is called the support of $X$ and is denoted by $S_X$. We can show this lemma by modifying the proof of \cite[Theorem 1.7.10]{Lali} or \cite[Theorem 3.16]{Steuding2} (see also \cite[Lemma 3.5]{Nakamura5}). 
\begin{lemma}
Let $g_n$ and $h_n$ be non-vanishing continuous bounded functions and $\{ X_n \}$ be a sequence of independent $H(D)$-valued random elements. Suppose that the series $\exp (\sum_{n=1}^{\infty} f_n (X_n)) - \exp ( \sum_{n=1}^{\infty} g_n (X_n))$ converges almost everywhere. Then the support of this sum is the closure of the set all $x \in H (D)$ which may be written as a convergent series 
\begin{equation}
x = \exp \Biggl(\sum_{n=1}^{\infty} f_n (x_n) \Biggr) - \exp \Biggl( \sum_{n=1}^{\infty} g_n (x_n) \Biggr) 
\mbox{ with }
x_n \in S_{X_n}. 
\label{eq:ao}
\end{equation}
\label{lem:aoyama}
\end{lemma}

\begin{proof}[Theorem \ref{mth:m1}]
Let $K$ be a compact subset of the critical strip $D$ with connected complement as introduced in Section 1. In view of Lemma \ref{lem:2}, there exist $\omega \in \Omega$ and a positive integer $N$ such that
$$
\sup_{s \in K} \biggl| \sum_{n>N}^{\infty} \log \biggl( 1 - \frac{\omega^j (n)}{p_n^s} \biggr) \biggr| < 
\frac{\varepsilon}{2} , \qquad 
\sup_{s \in K} \biggl| \sum_{n>N}^{\infty} \log \biggl( 1 - \frac{\omega^k (n)}{p_n^s} \biggr) \biggr| < 
\frac{\varepsilon}{2} .
$$
Suppose that $\omega (n) = 1$ for all $1 \le n \le N$. Then we have
\begin{equation*}
\begin{split}
&\sup_{s \in K} \Bigl| \log \zeta (s,\omega^j) - \log \zeta (s,\omega^k) \Bigr| = \sup_{s \in K} \Biggl| \sum_{n>N}^{\infty} \log \biggl( 1 - \frac{\omega^j (n)}{p_n^s} \biggr) -
\sum_{n>N}^{\infty} \log \biggl( 1 - \frac{\omega^k (n)}{p_n^s} \biggr) \Biggr| < \varepsilon .
\end{split}
\end{equation*}
By Lemma \ref{lem:aoyama}, the above inequality and $|e^{y+\varepsilon} - e^y| = |e^y| |e^{\varepsilon}-1|$, where $y \in {\mathbb{C}}$, the support of the random element $\zeta (s,\omega^j) - \zeta (s,\omega^k)$ contains a function $x(s) \in H(D)$ which satisfies $\sup_{s \in K} |x(s)| < \varepsilon$. We denote by $\Phi$ the set of functions $\phi \in H(D)$ such that $\sup_{s \in K} |\phi (s) - x(s)| < 2 \varepsilon$.  By using Limit theorem \ref{th:lim1}, the triangle inequality, and the property of support and the set $\Phi$ is open, we have
\begin{equation*}
\begin{split}
&\liminf_{T \rightarrow \infty} \nu_T \Bigl\{ \sup_{s \in K} \bigl| \zeta (s+ij\tau) - \zeta (s+ik\tau) \bigr| 
< \varepsilon \Bigr\} \\
\ge &\liminf_{T \rightarrow \infty} \nu_T \Bigl\{ \sup_{s \in K} \bigl| \zeta (s+ij\tau) - \zeta (s+ik\tau) - x(s) \bigr| < 2 \varepsilon \Bigr\}  
= P_T (\Phi) \ge P (\Phi)> 0 .
\end{split}
\end{equation*}
This inequality proves Theorem \ref{mth:m1}. 
\end{proof}

\section{In the region of absolute convergence}
In this section, we show that the generalized strong recurrence in the region of absolute convergence holds for any $d \in {\mathbb{R}}$. Let ${\mathcal{D}} := \{s \in {\mathbb{C}} : \Re (s) > 1\}$ and ${\mathcal{K}}$ be a compact subset contained in the strip ${\mathcal{D}}$. The following theorem  with $d=0$ should be compared with effective upper bounds for the almost periodicity of polynomial Euler products in the half-plane of absolute convergence \cite[Theorem 1]{Gir} (see also \cite[Theorem 9.6]{Steuding2}).

\begin{theorem}
For all $d \in {\mathbb{R}}$ and for any $\varepsilon > 0$ and ${\mathcal{K}}$, we have
\begin{equation}
\liminf_{T \rightarrow \infty} \nu_T \Bigl\{ \sup_{s \in {\mathcal{K}}} \bigl| \zeta (s+ i \tau) - \zeta (s+ id \tau) \bigr| < \varepsilon \Bigr\} > 0 .
\label{eq:1}
\end{equation}
\label{th:m1}
\end{theorem}

In order to prove Theorem \ref{th:m1}, we need the following lemma. 
\begin{lemma}{\rm{(see \cite[lemma 2.4]{Pan}).}}
Let ${\mathcal{P}}$ be the set of all primes. For arbitrary real irrational number $d$, there exists a finite set of primes $A_d$ containing at most two elements such that the set $\{ \log p_n \}_{{\mathcal{P}} \setminus A_d} \cup \{ \log p_n^d \}_{\mathcal{P}}$ is linearly independent over ${\mathbb{Q}}$. 
\label{lem:pan}
\end{lemma}

\begin{proof}[of Theorem \ref{th:m1}]
The idea of the proof is partly comes from the proof of \cite[Theorem 4]{Mishou}. By the definition of ${\mathcal{K}}$, for any $\tau$ and $d$, there exists a positive integer $N$ such that 
\begin{equation}
\sup_{s \in {\mathcal{K}}} \Biggl| \sum_{n \le N} \log(1-p_n^{-s-id\tau})^{-1} - \log \zeta (s+id\tau) \Biggr|
 < \varepsilon. 
\label{eq:zesyuu}
\end{equation}
We take a sufficiently small positive $\delta$ such that the inequality
\begin{equation}
\sup_{s \in {\mathcal{K}}} \Biggl| \sum_{n \le N} \log (1- p_n^{-s-i\tau})^{-1} - \sum_{n \le N} \log (1- p_n^{-s})^{-1} \Biggr| < \varepsilon. 
\label{eq:delta}
\end{equation}
holds when $\tau$ satisfies $|\exp(i\tau \log p_n) - 1| < \delta$. By Kronecker's approximation theorem, the set of $\tau$ which satisfies $|\exp(i\tau \log p_n) - 1| < \delta$ has a positive lower density. Therefore we have Theorem \ref{th:m1} when $d =0$ by (\ref{eq:zesyuu}). 

Next suppose $d = j/k$, where $j$ and $k$ are coprime integers. In this case, we have
\begin{equation}
\begin{split}
|\exp(i(j/k)\tau \log p_n) - 1| = & |\exp(i\tau \log p_n) - 1| 
\frac{|\sum_{m=0}^{j-1} \exp(i(m/j)\tau \log p_n)|}{|\sum_{m=0}^{k-1} \exp(i(m/k)\tau \log p_n)|} .
\label{eq:hyouka}
\end{split}
\end{equation}
Suppose $\tau$ satisfies $|\exp(i\tau \log p_n) - 1| < \min \{ \delta, \, |k/j|\delta \}$. By (\ref{eq:hyouka}), we have 
$$
\sup_{s \in {\mathcal{K}}} \Biggl| \sum_{n \le N} \log \bigl(1- p_n^{-s-i(j/k)\tau}\bigr)^{-1} - \sum_{n \le N} \log (1- p_n^{-s})^{-1} \Biggr| < \varepsilon .
$$
Therefore, by (\ref{eq:delta}) and the triangle inequality, we have
$$
\sup_{s \in {\mathcal{K}}} \Biggl| \sum_{n \le N} \log \bigl(1- p_n^{-s-i(j/k)\tau}\bigr)^{-1} - \sum_{n \le N} \log (1- p_n^{-s-i\tau})^{-1} \Biggr| < 2\varepsilon. 
$$
Thus we obtain Theorem \ref{th:m1} when $d = j/k$ by the above formula and (\ref{eq:zesyuu}). 

Finally, suppose $d$ is irrational. Put $A := \{ a_1, a_2 \}$ in Lemma \ref{lem:pan}, and $\log a_h := \sum_{n=1}^l\alpha_{h,n}\log p_n$, where $h=1,2$ and $\alpha_{h,n}$ is zero or an irreducible fraction. When $\tau$ satisfies $|\exp(i\tau \log p_n) - 1| < \delta$ for $p_n \in {\mathcal{P}} \setminus A$ and $|\exp(i\tau \log p_n^d) - 1| < \delta$ for $p_n \in {\mathcal{P}}$, we have
\begin{equation*}
\begin{split}
&|\exp(i\tau \log a_h) - 1| \\ = & \, \bigl| \exp(i\tau \log a_h) - \exp(i\tau {\textstyle{\sum_{n=1}^{l-1}}} \alpha_{h,n} \log p_n) + \exp(i\tau {\textstyle{\sum_{n=1}^{l-1}}} \alpha_{h,n} \log p_n) -1 \bigl| \\ \le & \,
| \exp(i\tau \alpha_{h,l} \log p_l) - 1 | + 
\bigl| \exp(i\tau {\textstyle{\sum_{n=1}^{l-1}}} \alpha_{h,n} \log p_n) -1 \bigr| 
\le \delta {\textstyle{\sum_{n=1}^l}} |\alpha_{h,n}|.
\end{split}
\end{equation*}
by the manner similar to the proof of (\ref{eq:hyouka}). Hence we have
$$
\sup_{s \in {\mathcal{K}}} \Biggl| \sum_{n \le N} \log \bigl(1- p_n^{-s-id\tau}\bigr)^{-1} - \sum_{n \le N} \log (1- p_n^{-s-i\tau})^{-1} \Biggr| < 2\varepsilon. 
$$
by the method similar to the case $d = j/k$. Therefore we obtain Theorem \ref{th:m1} when $d$ is irrational. 
\end{proof}

\subsection*{Note}
This paper was submitted to a journal on 12 March 2010. Corollary 1.2 of this paper coincides with Corollary 2 of R.~Garunkstis, Self-approximation of Dirichlet $L$-functions, arXiv:1006.1507, for $\zeta (s)$. However the proof in this paper is Bagchi's method.



\begin{thebibliography}{1}
  \bibitem{BagchiZ}
{\sc B.~Bagchi}, `A joint universality theorem for Dirichlet $L$-functions', {\em Math.~Z.} 181 (1982), no.~3, 319--334.
  \bibitem{Denjoy}
{\sc A.~Denjoy}, `L'Hypoth\'{e}se de Riemann sur la distribution des z\'{e}ros de $\zeta (s)$, reli\'{e}e
\`{a} la th\'{e}orie des probabilit\'{e}s', {\em Comptes Rendus Acad. Sci. Paris} 192 (1931), 656--658.
  \bibitem{Edwards}
{\sc H.~M.~Edwards}, {\em Riemann's zeta function}, Pure and Applied Mathematics, Vol. 58. Academic Press, New York-London, 1974. 
  \bibitem{Gir}
{\sc E.~Girondo} and {\bibname J.~Steuding}, `Effective estimates for the distribution of values of Euler products', {\em Monatsh.~Math.} 145 (2005), no.~2, 97--106. 
  \bibitem{Lali}
{\sc A.~Laurin\v cikas}, 
{\em Limit Theorems for the Riemann Zeta-function}, Kluwer Academic Publishers, 1996.
  \bibitem{Mishou}
{\sc H.~Mishou}, `The joint value distribution of the Riemann zeta function and Hurwitz zeta functions.~II', {\em Arch.~Math.~(Basel)} 90 (2008), no.~3, 230--238. 
  \bibitem{Nakamura1} 
{\sc T.~Nakamura}, `The joint universality and the generalized strong recurrence for Dirichlet {\it{L}}-functions', {\em Acta Arith.} 138 (2009), no.~4 357--362. 
  \bibitem{Nakamura2} 
{\sc T.~Nakamura}, `Some topics related to universality of $L$-functions with an Euler product', to appear in {\em Analysis}. 
  \bibitem{Nakamura5}
{\sc T.~Nakamura}, `Value distribution of zeta functions associated to symmetric matrices', {\em preprint.}
  \bibitem{Pan}
{\sc \L.~Pa\'{n}kowski}, `Some remarks on the generalized strong recurrence for $L$-functions', {\em New Directions in Value Distribution Theory of zeta and L-Functions}, Proceedings of W\"{u}rzburg Conference, October 6-10, 2008, Shaker Verlag, (2009), 305-315.. 
  \bibitem{Steuding2}
{\sc J.~Steuding}, {\em Value Distributions of $L$-functions}, Lecture Notes in Mathematics 1877, Springer-Verlag, 2007.
  \bibitem{Tit}   
{\sc E.~C.~Titchmarsh}, {\em The theory of the Riemann zeta-function}, Second edition. Edited and with a preface by D.~R.~Heath-Brown. The Clarendon Press, (Oxford University Press, New York, 1986). 
  \bibitem{Voronin}
{\sc S.~M.~Voronin}, {\em Theorem on the universality of the Riemann zeta-function}, Izv.~Akad.~Nauk.~SSSR Ser.~Mat. 39 (1975), 475-486 (in Russian); Math.~USSR Izv. 9 (1975), 443-453. 
  \end{thebibliography}
\end{document}